\documentclass{article}
\usepackage{amssymb,amsmath}
\usepackage{graphics}
    
\def\qed{\hfill {\hbox{${\vcenter{\vbox{               
   \hrule height 0.4pt\hbox{\vrule width 0.4pt height 6pt
   \kern5pt\vrule width 0.4pt}\hrule height 0.4pt}}}$}}}
   
\def\tr{\triangleright}
\def\tl{\triangleleft}

{\qed\par\medskip}

\newtheorem{definition}{Definition}

\newtheorem{example}{Example}

\author{
{\begin{tabular}{c} Richard Henderson$^{\ast}$ \\
\small{\texttt{rth@redhat.com}}\end{tabular}}
\and
{\begin{tabular}{c} Todd Macedo$^{\dagger}$ \\
\small{\texttt{tmace001@student.ucr.edu}}\end{tabular}}
\and 
{\begin{tabular}{c} Sam Nelson$^{\dagger}$ \\
\small{\texttt{knots@esotericka.org}}\end{tabular}}
\and
\small{${\ }^{\ast}$ Red Hat, Inc., 444 Castro Street, Suite 1200, Mountain View, CA, 94041}
\and
\small{${\ }^{\dagger}$ University of California, Riverside, 900 University Avenue, Riverside, CA, 92521 }
}
 
\date{}
 
\title{\Large \textbf{Symbolic computation with finite quandles}}
 
\begin{document}

\maketitle

\begin{abstract}
\noindent Algorithms are described and \textit{Maple} implementations are 
provided for 
finding all quandles of order $n$, as well as computing all homomorphisms 
between two finite quandles or from a finitely presented quandle (e.g., a knot 
quandle) to a finite quandle, computing the automorphism group of a finite 
quandle, etc. Several of these programs work for arbitrary binary
operation tables and hence algebraic structures other than quandles. We also 
include a stand-alone C program which finds quandles of order $n$ and provide
URLs for files containing the results for $n=6,$ 7 and 8. 
\end{abstract}



\section{Introduction}

In 1980, David Joyce introduced a new algebraic structure dubbed the
\textit{quandle}. Quandles are tailor-made for defining invariants of knots
since the quandle axioms are essentially the Reidemeister moves written
in algebra. Associated to any knot diagram, there is a quandle called
the \textit{knot quandle} which is a complete invariant of knot type up
to homeomorphism of topological pairs.

The history of quandle theory is a story of rediscovery and reinvention.
Quandles and their generalization, racks, have been independently invented and
studied by numerous authors (\cite{B}, \cite{FR}, \cite{J}, \cite{M}, etc.)
and classification results for various subcategories of quandles have
been obtained by various authors (\cite{G}, \cite{N}.) In \cite{HN}, the third 
listed author and a coauthor described a way of representing finite quandles 
as matrices and implemented algorithms for finding all finite quandles, 
removing isomorphic quandles from the list, and computing the automorphism 
group of each quandle. As we later learned, some of our work has duplicated 
the efforts of others (\cite{R}, \cite{LR}, \cite{C}.)

This paper is intended to reduce future duplication of effort by describing
the algorithms for computation with finite quandles implemented in \cite{HN} 
and other recent projects, as well as an improved algorithm for finding 
quandle matrices. The C source for our implementation of this algorithm as 
well as \textit{Maple} implementations of algorithms for computing with 
finite quandles 
and the lists of quandle matrices of order 6, 7 and 8 are available for 
download at \texttt{http://www.esotericka.org/quandles}. Additional maple code 
corresponding to current and future projects will be made available at the 
same site, such as an algorithm for finding all Alexander presentations of 
a quandle when such exist \cite{MNT}.

\section{Quandles, quandle matrices, and homomorphisms} \label{sec2}

\begin{definition}\textup{
A \textit{quandle} is a set $Q$ with a binary operation $\tr:Q\times Q\to Q$
satisfying
\newcounter{qax}
\begin{list}{(\roman{qax})}{\usecounter{qax}}
\item{for every $x\in Q$ we have $x\tr x=x$,}
\item{for every $x,y\in Q$ there is a unique $z\in Q$ such that $x=z\tr y$, 
and}
\item{for every $x,y,z\in Q$ we have $(x\tr y) \tr z = (x\tr z)\tr (y\tr z)$.}
\end{list}
If $(Q,\tr)$ satisfies (ii) and (iii), $Q$ is a \textit{rack}.
}
\end{definition}

Axiom (ii) says that $\tr$ is right-invertible; for every $y\in Q$, the map
$f_y:Q\to Q$ defined by $f_y(x)=x\tr y$ is a bijection (indeed, a quandle
automorphism). Denote the inverse map $f_y^{-1}(x)=x\tl y$. Then $(Q,\tl)$
is also a quandle, called the \textit{dual} of $(Q,\tr)$; not only is 
$\tl$ self-distributive, but it is an easy exercise to check that $\tr$ and
$\tl$ distribute over each other.

Standard examples of quandles include groups, which are quandles under 
conjugation $g\tr h=h^{-1}gh$ as well as $n$-fold conjugation 
$g\tr h=h^{-n}gh^n$, denoted $\mathrm{Conj}(G)$ and $\mathrm{Conj}_n(G)$
respectively, and Alexander quandles, which are modules over the ring
$\Lambda=\mathbb{Z}[t^{\pm 1}]$ of Laurent polynomials in one variable with 
integer coefficients, with quandle operation given by
\[ x\tr y = tx + (1-t)y.\]

A finite quandle $Q$ may be specified by giving its \textit{quandle matrix}
$M_Q$, which is the matrix obtained from the operation table of $Q=\{x_1,x_2,
\dots,x_n\}$ (where the entry in row $i$ column $j$ is $x_i\tr x_j$) by 
dropping the $x$s and keeping only the subscripts. In \cite{HN} it is noted 
that, unlike arbitrary binary operation tables or indeed even rack tables, 
quandle axiom (i) permits us to deduce the column and row labels from the 
elements along the diagonal of a quandle matrix.

\begin{example} \textup{
Let $Q=R_4$, the dihedral quandle of order 4, which has underlying set 
$Q=\{x_1=0,x_2=1,x_3=2,x_4=3\}$ with quandle operation $x_i\tr x_j= x_{2j-i \ 
(\mathrm{mod} 4)}.$ Then $Q$ has operation table
\[
\begin{array}{c|cccc}
    & x_1 & x_2 & x_3 & x_4 \\ \hline
x_1 & x_1 & x_3 & x_1 & x_3 \\
x_2 & x_4 & x_2 & x_4 & x_2 \\
x_3 & x_3 & x_1 & x_3 & x_1 \\
x_4 & x_2 & x_4 & x_2 & x_4 
\end{array}
\quad \mathrm{and \ hence \ matrix} \quad
M_{R_4}=
\left[
\begin{array}{rrrr}
1 & 3 & 1 & 3 \\
4 & 2 & 4 & 2 \\
3 & 1 & 3 & 1 \\
2 & 4 & 2 & 4
\end{array}
\right].
\] }
\end{example}

A map $\phi:Q\to Q'$ from a quandle $Q=\{x_1,\dots,x_n\}$ to a 
quandle $Q'=\{y_1,\dots, y_m\}$ may be represented by a vector
$v=(\phi(x_1),\phi(x_2),\dots, \phi(x_n))\in {Q'}^{n}$. Such a vector
$v$ then represents a homomorphism iff $\phi(x_i\tr x_j)=\phi(x_i)\tr 
\phi(x_j)$, that is, iff we have
\[v[A[i,j]]=B[v[i],v[j]]\]
for all $x_i,x_j\in Q$ where $A=M_Q$, $B=M_{Q'}$, and the notation
$M[i,j]$ indicates the entry of $M$ in row $i$ column $j$.

In \cite{FR}, presentations of quandles by generators and relations are
defined. In \cite{N2}, it is observed that all finitely presented quandles
may be written with a \textit{short form} presentation in which every relation
is of the form $a=b\diamond c$ where $\diamond \in \{\tr,\tl\}$. In 
particular, a knot quandle has a presentation with $n$ such short relations
where $n$ is the number of crossings in the diagram. Moreover, we may assume 
(rewriting if necessary) that every relation is written in the form
$a=b\tr c$ and that no two relations of the form $a=b\tr c$ and $a'=b\tr c$ 
are present, since if $a=b\tr c$ and $a'=b\tr c$ are both present we can
replace every instance of $a'$ with $a$ and remove the generator $a'$
without changing the presented quandle; in particular, if our quandle is a 
knot quandle, Reidemeister type I moves\footnote{Reidemeister 
moves are described in \cite{N2} and many other works.} induce such a 
replacement. 

\begin{definition}\textup{
Let $Q=\langle 1,2,\dots, n \ | \ a_1=b_1\tr c_1, \dots, a_m=b_m\tr c_m , 
m\le n^2\rangle$ be a short form quandle presentation such that no two 
relations 
of the form $a_i=b_i\tr c_i$ and $a_j=b_i\tr c_i$ with $a_i\ne a_j$
are present. The matrix $MP\in M_n(\mathbb{Z})$ with
\[
MP[i,j]=\left\{
\begin{array}{ll}
k & \quad k = i\tr j  \mathrm{ \ a \ listed \ relation} \\
0 & \quad \mathrm{otherwise}
\end{array}
\right.
\]
is the \textit{matrix of the presentation $Q$}. Note that a quandle matrix
for a finite quandle is the matrix of a presentation of a finite quandle,
so this definition generalizes the notion of quandle matrices to finitely 
presentable quandles.}
\end{definition}

\begin{example}\textup{
\[\raisebox{-0.5in}{\includegraphics{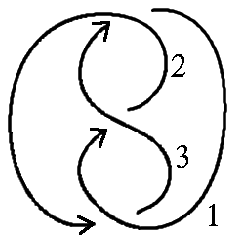}} \qquad 
MP_{KQ} = \left[
\begin{array}{rrr}
0 & 3 & 0 \\
0 & 0 & 2 \\
1 & 0 & 0 
\end{array}
\right]\]
The pictured trefoil knot diagram has quandle presentation 
$\langle 1,2,3 \ |  1=2\tr 3, 2=3\tr 1, 3=1\tr 2\rangle.$
The relations are determined at a crossing by looking in the positive 
direction of the overcrossing strand indicated by the given orientation; 
the relation is 
\begin{center}
(left-hand undercrossing) = (right-hand undercrossing) $\tr$ (overcrossing).
\end{center}
See \cite{FR} or \cite{N2} for more.  }
\end{example}

This matrix representation gives us a convenient way to do computations
involving quandles, including the quandle counting invariant for knot 
quandles or other short form quandles with respect to a finite target quandle. 
The next section describes algorithms for doing computations with quandles 
and refers to implementations in \textit{Maple} \cite{W} and C \cite{H}.

\section{Algorithms} \label{al}

The goal of the computations in \cite{HN} was to find all quandles of 
a given finite order. Originally, we wrote separate programs for each value of
$n$; \cite{W} includes one example of such an implementation,
\texttt{quandleslist5}. We later wrote a more general program which works 
for arbitrary $n$, though due to the large number of columns to be checked,
for values of $n\ge 6$ we decided to implement a stand-alone version suitable 
for distributed computing. 

The algorithm implemented in \texttt{quandleslist} takes a number $n$
and generates a list of all $n\times n$ standard form quandle matrices. 
A matrix $M\in M_n (\mathbb{Z})$ is a quandle matrix in standard form iff it
satisfies the following three conditions:
\newcounter{ql}
\begin{list}{(\roman{ql})}{\usecounter{ql}}
\item{for $i\in \{1,\dots, n\}$, $M[i,i]=i$,}
\item{every column in $M$ is a permutation of $\{1,\dots, n\}$, and}
\item{for every triple $1\le i,j,k \le n$ we have
$M[M[i,j],k]=M[M[i,k],M[j,k]]$.}
\end{list}

To guarantee that conditions (i) and (ii) are satisfied, we start by 
getting a list of all permutations of $\{1,\dots, n\}$. The program
\texttt{listperms} takes a number $n$ and produces a list of all permutations
$\rho\in \Sigma_n$, represented as vectors $[\rho(1),\rho(2),\dots, \rho(n)]$,
in the dictionary order.

The $i$th column in a standard form quandle matrix has entry $i$ in 
the $i$th position. The program \texttt{listpermsi} takes a pair of positive
integers $(n,i)$ and outputs a list of all permutations of $\{1,2,\dots, n\}$
$\rho\in \Sigma_n$ which fix the element $i$.

To test quandle axiom (iii), we note that the first time any triple $(i,j,k)$ 
fails to satisfy the axiom, we can exit the program and report that the matrix 
is not a quandle. This is implemented in \texttt{q3test}.

For a fixed value of $n$, we can then simply run over a series of nested
loops, testing each resulting matrix for quandle axiom (iii), since
by construction axioms (i) and (ii) are already satisfied. The program 
\texttt{quandleslist5} is an example of this.

The program \texttt{quandleslist} finds a list of all quandle matrices of a 
given size $n$. To find all $n\times n$ quandle matrices for arbitrary $n$, 
\texttt{quandleslist} finds all control vectors $v[i]$ with $n$ entries using 
\texttt{listmaps}, a program which takes two inputs $a$ and $b$ and outputs a 
list of all $a$-component vectors with entries in $\{1,\dots, b \}$. Each 
entry in the control vector corresponds to a column in the output matrix; for 
each such control vector, an $n\times n$ matrix $M[i,j]$ is produced whose 
$i$th column is $L[n,i][v[i]]$, where $L[n,i]$ is the
output of \texttt{listpermsi(n,i)}. These matrices are then tested for quandle
axiom (iii) using \texttt{q3test}. For completeness, we include a program 
which tests a matrix for all three quandle axioms, \texttt{qtest}.

Since every $n$-component vector with entries in $\{1,\dots, m\}$ can be
interpreted as a map from $\{1,\dots, n\}$ to $\{1,\dots, m\}$, we can use
\texttt{listmaps} to compute the set of all homomorphisms
from one finite quandle to another. Let $A\in M_n(\mathbb{Z})$ be an
$n\times n$ quandle matrix and $B\in M_m(\mathbb{Z})$ an $m\times m$ quandle 
matrix. Then the vector $v\in \mathbb{Z}^n$, $1\le v[i]\le m$ represents
a quandle homomorphism $v:A\to B$ iff 
\[ v[A[i,j]]=v(i\tr j) =v(i)\tr v(j)=B[v[i],v[j]],\]
as noted in section \ref{sec2}. The program \texttt{homtest} takes two 
quandle matrices and a vector and reports whether the vector represents a 
quandle homomorphism or not. 

The program \texttt{homtest} handles the case where $A$ is 
either a finite quandle matrix or a presentation matrix for a finitely 
presented quandle; in the former
case, the program simply tests whether the assignment of generators 
$\{1,\dots, n\}$ in the quandle with presentation matrix $A$ to elements
$\{1,\dots, m\}$ in the finite quandle $B$ satisfy the relations defining
$A$ by ignoring any zero entries in $A$. 

We make use of \texttt{nextmap}, a procedure which takes as input a vector $v$
and number $n$ and returns the next $m$-component vector with entries in
$\{1,2,\dots, m\}$ in the dictionary order, to get a list of all homomorphisms 
from
the quandle with matrix $A$ to the quandle with matrix $B$ in the program
\texttt{homlist}. The program \texttt{homcount} counts the number of 
homomorphisms from one finite quandle to another. If $A$ is a knot quandle
presentation matrix, then \texttt{homcount} computes the quandle counting
invariant, i.e., the number of quandle colorings of the knot diagram 
defining $A$ by the finite quandle $B$. Alternate methods of computing the
quandle counting invariant for finite Alexander quandles are described in 
\cite{DL}.

After the first version of this paper was completed, we implemented a much 
faster algorithm for finding quandle homomorphisms, \texttt{homlist2}. This
program uses a $|B|$-component vector with entries in $\{0,1,\dots,|A|\}$
as a template for a homomorphism, with 0 entries acting as blanks to be
filled in. The program keeps a working list of such templates, systematically
filling in zero values and propagating the value through the template using
\texttt{homfill}. The procedure \texttt{homfill} takes as input a quandle 
matrix $B$, a quandle presentation matrix $A$ and a template vector $v$
and systematically checks every pair of entries for the quandle homomorphism
condition $v[A[i,j]]=B[v[i],v[j]]$, filling in zeroes where possible and 
reporting ``false'' if a contradiction is found. 

Since an isomorphism is a bijective homomorphism, and a bijective map
is represented by a permutation $v:\{1,\dots, n\}\to \{1,\dots,n\}$,
we can test whether two quandles given by matrices are isomorphic
by running through the list of permutations of order $n$
and testing to see whether any are homomorphisms. The program 
\texttt{isotest} returns ``true'' the first time it finds an isomorphism
and ``false'' if it gets through all $n!$ permutations without finding
an isomorphism.\footnote{A faster version of this program using orbit
decompositions of finite quandles is described in \cite{NW}.}

Replacing \texttt{listmaps} in \texttt{homlist2} with \texttt{permute(n)}
and setting $B=A$ gives us the automorphism group of the quandle with 
matrix $A$, \texttt{autlist}, represented as a list of permutation vectors.

In \cite{HN}, a slightly different method of determining the automorphism 
group of a quandle was used. Specifically, permuting the entries of
a quandle matrix $A$ by a permutation $\rho$ applies an isomorphism to the
defined quandle, but the new matrix now has its rows and columns out of 
order. To restore the order, we conjugate by the matrix of the permutation; 
the resulting matrix was called $\rho(A)$ in \cite{HN}. In particular,
a permutation $\rho$ is an automorphism of $A$ iff $\rho(A)=A$. To compute
$\mathrm{Aut}(A)$ in \cite{HN}, we ran a loop over the list of permutations 
given by \texttt{listperms} and noted which ones preserved the original 
matrix $A$. Here, we include a program \texttt{stndiso} which computes the 
standard form matrix $\rho(A)$ given a quandle matrix $A$ and a vector $v$ 
representing the permutation $\rho$.

Finally, once we have a list of quandle matrices of order $n$, we want
to choose a single representative from each isomorphism class. The program 
\texttt{reducelist} takes a list of quandle matrices and compares them
pairwise with \texttt{isotest}, removing isomorphic copies and outputting a 
trimmed list. The program \texttt{reducelist} works for short lists; an 
improved algorithm, implemented as \texttt{reducelist2}, is better for 
longer lists, but neither is sufficient to reduce the rather lengthy lists 
of quandles of order 7 and 8 in a reasonable amount of time.

We note that several of these programs, notably \texttt{homtest},
\texttt{homlist}, \texttt{homlist2}, \texttt{homcount}, \texttt{isotest}, 
\texttt{autlist}, and \texttt{reducelist} are not quandle-specific but apply 
as written to any binary operation defined using a matrix as operation table. 
These facts are exploited in \cite{MNT}, in which the authors give a program 
which determines all Alexander structures on a quandle, if there are any, using
matrices to represent the Cayley table of an abelian group.

We have also implemented a stand-alone version of \texttt{quandleslist}, 
written in C (see \cite{H}); it writes a list of quandle matrices in 
\textit{Maple} format to an output file.

In our initial version of the stand-alone program, several 
instances of the program could be run in parallel on networked machines 
using a control file to ensure that separate instances do not repeat the 
same computations. However, sufficient improvements were made to the 
algorithm by pruning the search space that the current version can handle 
the $n=8$ case on a single processor, though the $n=9$ case is still out of 
reach even with a large network.  

The first improvement was to introduce a partial test versus axiom (iii)
after generation of each column.  In many cases we can find entries that
violate the axiom well before the entire matrix is generated, which
allows vast portions of the search space to be pruned.

The second improvement was to notice when all of the interior coordinate
values as well as the left-hand side value of the axiom (iii) equality
have been computed, but the right-hand side value has not.  In this case
we can constrain a row of a future column to be equal to the left-hand
side value.  This reduces the number of rows that must be permuted when
searching that column, which further prunes the search space.  The
pruning effect is magnified the earlier these constraints are added.
For example, with n=7 and all else held equal, adding a single constraint 
to column 3 saves  $(6! -5!)*((6!)^4)$ or $2.3\times 10^{16}$ tests, whereas 
adding a constraint to column 7 saves only $6! - 5!$ or 600 tests.

The effect of the two improvements can be seen in the following tables. It 
is interesting to note that although there is nothing in the program
to prevent it (and reasonable amount of code to encourage it), we never
add constraints to column 2, nor do we ever add more than one constraint
per column, nor do we detect addition of conflicting constraints.

\begin{center}

\begin{tabular}{|l|llll|} \hline
$n$ &                    5  &  6     &      7 & 8\\ \hline
total search space & $8.0\times 10^6$ & $3.0\times 10^{11}$ 
& $1.0\times 10^{20}$ & $4.1\times 10^{30}$ \\
with early testing       & 8400  & 715680 & $1.0\times 10^8$  & n/a \\
with forward propagation & 1154  & 53500  & $5.0\times 10^6$ 
& $7.7\times 10^8$ \\
total quandles           & 404   & 6658   & 152900 & 5225916 \\ \hline
\end{tabular}

\smallskip

Table 1: Number of complete matrices tested.

\medskip

\begin{tabular}{|l|lllll|} \hline
column         & 2      &  3       & 4       & 5       &   6 \\ \hline
$n=5  $          & 186    &  8736    & 14626   &  -       & - \\
$n=6 $           & 4728   &  1090404 & 8418374 & 1187556 & - \\
$n=7$            & 154680 &  $2.3\times 10^8 $  & $1.8\times 10^{10}$  & 
$3.8\times 10^9$   & $1.8\times 10^8$ \\ \hline
\end{tabular} 

\smallskip

Table 2: Number of tests pruned with early testing. 

\medskip

\begin{tabular}{|l|lllll|} \hline
column &  3       &  4      & 5     & 6     & 7 \\ \hline
$n=5$  &  164     &  179    & 290   & -     & - \\
$n=6$  &  3558    &  4396   & 3348  & 5020  & - \\
$n=7$  &  115872  &  228384 & 91452 & 82910 & 117430 \\ \hline
\end{tabular}

\smallskip

Table 3: Number of times columns constrained beyond axiom (i).

\end{center}


\end{document}